\documentclass[12pt]{article}
	\usepackage{amsmath,fullpage,amsthm,amssymb,xcolor,graphicx}
\newtheorem{prop}{Proposition}[section] 
	\theoremstyle{theorem}
	\newtheorem{thm}{Theorem}[section]
	\theoremstyle{corollary}
	\newtheorem{cor}{Corollary}[section]
	\theoremstyle{lemma}
	
	\theoremstyle{definition}
	\newtheorem{defn}{Definition}[section]
	\theoremstyle{proof}
	\theoremstyle{remark}
	\newtheorem{remark}{Remark}[section]
	\theoremstyle{example}
	\newtheorem{exa}{Example}[section]
	\theoremstyle{observation}

	\title{On Weaving Generalized Frames and Generalized Riesz Bases}
\author{Deepshikha \thanks{Department of Mathematics, Shyampur Siddheswari Mahavidyalaya,University of Calcutta-711312, India. Email: dpmmehra@gmail.com}\  \and Aniruddha Samanta \thanks{Department of Mathematics, Indian Institute of Technology Kharagpur, Kharagpur 721302, India. Email: aniruddha.sam@gmail.com, aniruddha.samanta@iitkgp.ac.in }
}
\date{\today}
\begin{document}
\maketitle
\baselineskip=0.25in

\begin{abstract}
Weaving frames have potential application in wireless sensor networks that require distributed processing of signal under different frames. In this paper, we study some new properties of weaving generalized frames (or  $g$-frames) and weaving generalized orthonormal bases (or $g$-orthonormal bases). It is shown that a $g$-frame and its dual $g$-frame are woven. Inter relation of optimal $g$-frame bounds and optimal universal $g$-frame bounds is studied. 
 Further, we present a characterization of weaving $g$-frames. Illustrations are given to show the difference in properties of weaving generalized Riesz bases and weaving Riesz bases.
\end{abstract}



\section{Introduction}
  The concept of frames was introduced by Duffin and Schaeffer in \cite{DS} and popularized by Daubechies, Grossman, and Meyer in \cite{DGM} when they showed the importance of frames in data processing. Frames can view as the generalization of bases but allow for over completeness. This redundancy of frames makes them play a vital role in numerous areas viz., noise reduction, sparse representations, image compression, signal transmission and processing, image processing, and wavelet analysis etc. Frames also help to spread the information over a wider range of vectors and thus it provide resilience against losses or noises. For basic theory and applications of frames, we refer  \cite{CK, OC2,D,FT,HL}. 
  
  The concept of  generalized frames or $g$-frames was introduced by Sun in \cite{S}. $G$-frames are generalization of many frames like ordinary frames, frames of subspaces, pesudo frames etc. These frames are useful in many applications. In this paper, we study weaving generalized frames and weaving generalized Riesz bases. The notion of weaving frames was introduced by Bemrose et al. in  \cite{BCGLL}. Weaving frames have potential applications in wireless sensor networks that require distributed signal processing under different frames, as well as preprocessing of signals using Gabor frames.
  
  \section{Preliminaries} 
  Throughout the paper, $\mathcal{H}$ is a separable Hilbert space and $\{\mathcal{H}_m\}_{m \in \mathbb{N}}$ is a sequence of subspaces of a separable Hilbert space. $\mathcal{L}(\mathcal{H}, \mathcal{H}_m)$ is the space of all linear bounded operators from $\mathcal{H}$ to $\mathcal{H}_m$. If $U \in \mathcal{L}(\mathcal{H}, \mathcal{K})$  then $U^*$ denotes the Hilbert-adjoint operator of $U$, where $\mathcal{H}$ and $\mathcal{K}$ are Hilbert spaces. w.r.t. is the abbreviation used for \emph{with respect to}.
  
  \subsection{Frames}
  Suppose $\{h_m\}_{m\in \mathbb{N}}$ is a countable sequence of vectors in $\mathcal{H}$. Then  $\{h_m\}_{m\in \mathbb{N}}$  is called a \emph{frame } (or \emph{ordinary frame}) for $\mathcal{H}$ if there exist positive constants $A\leq B$ such that for any $h\in\mathcal{H}$,
  \begin{align}\label{1.1}
  	A \|h\|^2\leq  \sum_{m \in \mathbb{N}} |\langle h, h_m\rangle|^2 \leq B \|h\|^2.
  \end{align}
  If $\{h_m\}_{m\in \mathbb{N}}$ satisfies only upper inequality in \eqref{1.1}, then it is called a \emph{Bessel sequence}
   and $B$ is called a \emph{Bessel bound}. 
  
  If a frame $\{h_m\}_{m\in \mathbb{N}}$ ceases to be a frame when an arbitrary element is removed then $\{h_m\}_{m\in \mathbb{N}}$ is called an \emph{exact frame}.
  
  Associated with a Bessel sequence $\{h_m\}_{m\in \mathbb{N}}$, the \emph{frame operator} $S:\mathcal{H}\rightarrow\mathcal{H}$ is defined by
  \begin{align*}
  	S(h)=\sum_{m \in \mathbb{N}}\langle h, h_m \rangle h_m .
  \end{align*}
  The frame operator $S$ is linear, bounded and self adjoint. If  the Bessel sequence $\{h_m\}_{m\in \mathbb{N}}$ is a frame then the frame operator is invertible. Using the frame operator, we have a series representation of each vector $h\in\mathcal{H}$ in terms of frame elements which is given by
  \begin{align*}
  	&h = \sum_{m \in \mathbb{N}} \langle h, S^{-1}h_m \rangle h_m
  	\ =\sum_{m \in \mathbb{N}}\langle h,h_m\rangle S^{-1}h_m .
  \end{align*}
  
  If a frame $\{h_m\}_{m\in \mathbb{N}}$ is not exact then there exist a frame $\{g_m\}_{m\in \mathbb{N}}$ other than $\{S^{-1}h_m\}_{m\in \mathbb{N}}$ such that 
  \begin{align*}
  	&h = \sum_{m \in \mathbb{N}} \langle h, g_m \rangle h_m,\ \forall h\in\mathcal{H}.
  \end{align*}
  Here $\{g_m\}_{m\in \mathbb{N}}$ is called a \emph{dual frame} of $\{h_m\}_{m\in \mathbb{N}}$. Thus a frame can provide more than one series representation of a vector in terms of the frame elements.
  
  A sequence $\{h_m\}_{m\in \mathbb{N}}\subset\mathcal{H}$ is called a \emph{Riesz basis} for $\mathcal{H}$ if $\{h_m\}_{m\in \mathbb{N}}$ is complete in $\mathcal{H}$, and there exist positive constants $A\leq B$ such that for any finite scalar sequence $\{c_m\}$, 
  \begin{align*}
  	A\sum |c_m|^2\leq\left\|\sum c_m h_m\right\|^2\leq B\sum |c_m|^2.
  \end{align*}
  Riesz bases are the images of orthonormal bases under bounded invertible opertaors  \cite{OC2}. Thus these can be viewed as generalization of orthonormal bases.

  \subsection{Weaving frames} 
  We start with the definition of weaving frames which was given by Bemrose et al. in \cite{BCGLL}. 
  
  \begin{defn}
  	Frames $\{\phi_{m}\}_{m\in \mathbb{N}}$ and $\{\psi_{m}\}_{m\in \mathbb{N}}$ for $\mathcal{H}$ are called \emph{woven} if there exist positive constants $A\leq B$ such that for any $\sigma\subseteq\mathbb{N}$, $\{\phi_{m}\}_{m\in \sigma}\cup\{\psi_{m}\}_{m\in \sigma^c}$ is a frame for $\mathcal{H}$ with lower frame bound $A$  and upper frame bound $B$. Each $\{\phi_{m}\}_{m\in \sigma}\cup\{\psi_{m}\}_{m\in \sigma^c}$ is called a \emph{weaving}.
  \end{defn}
  Bemrose et al. presented one interesting result in  \cite{BCGLL} which says that a Riesz basis and a frame (which is not a Riesz basis) are not woven.
  \begin{thm}\label{th-1c}\cite{BCGLL}
  Suppose $\{\phi_{m}\}_{m\in \mathbb{N}}$ is a Riesz basis and $\{\psi_{m}\}_{m\in \mathbb{N}}$ is a frame for $\mathcal{H}$. If $\{\phi_{m}\}_{m\in \mathbb{N}}$ and $\{\psi_{m}\}_{m\in \mathbb{N}}$ are woven, then $\{\psi_{m}\}_{m\in \mathbb{N}}$ is a Riesz basis.
  \end{thm}
  Following result presented in \cite{BCGLL} says that if two Riesz bases are woven, then every weaving is a Riesz basis.
  \begin{thm}\label{th-2c}\cite{BCGLL}
  Suppose  $\{\phi_{m}\}_{m\in \mathbb{N}}$ and $\{\psi_{m}\}_{m\in \mathbb{N}}$ are Riesz bases for $\mathcal{H}$ and there is a uniform constant $A>0$ so that for any $\sigma\subset\mathbb{N}$, $\{\phi_{m}\}_{m\in \sigma}\cup\{\psi_{m}\}_{m\in \sigma^c}$ is a frame with lower frame bound $A$. Then for any $\sigma\subset\mathbb{N}$, $\{\phi_{m}\}_{m\in \sigma}\cup\{\psi_{m}\}_{m\in \sigma^c}$ is a Riesz basis.
  \end{thm}
   Many interesting properties of weaving frames were studied by Casazza et al. in \cite{CL}. Then the notion of weaving frames in Hilbert spaces was extended to Banach spaces in \cite{CFL}. A characterization for the weaving of approximate Schauder frames in terms of $C$-approximate Schauder frame was presented. The concept of weaving frames in different settings were studied by many authors in \cite{DV, DV2, KB,  VD1, VD2}.

  \subsection{Generalized frames in Hilbert spaces}
  Sun \cite{S} gave the concept of generalized frames or $g$-frames which is the generalization of ordinary frames, fusion frames, bounded quasi-projectors etc., see \cite{CK1,FM, KM,SO}. 
  \begin{defn} \cite{S}
  	A sequence  $\{\Lambda_m\in\mathcal{L}(\mathcal{H},\mathcal{H}_m):m \in \mathbb{N}\}$ is a  \emph{generalized frame} (or \emph{$g$-frame}) for $\mathcal{H}$ w.r.t. $\{\mathcal{H}_m\}_{m\in\mathbb{N}}$  if there exist  positive constants $A\leq B$ such that
  	\begin{align} \label{eq2.2}
  		A\|h\|^2 \leq \sum\limits_{m\in\mathbb{N}}\|\Lambda_m h\|^2 \leq B\|h\|^2,\ \forall h\in\mathcal{H}.
  	\end{align}
  \end{defn}
  The constants $A$ and $B$ are called lower and upper $g$-frame bounds, respectively. The supremum of all lower $g$-frame bounds is called the \emph{optimal lower $g$-frame bound}, and  the infimum of all upper $g$-frame bounds is called the \emph{optimal upper $g$-frame bound}.
   
  If $\{\Lambda_m\}_{m\in\mathbb{N}}$ satisfies the upper inequality in \eqref{eq2.2} then it is called a \emph{$g$-Bessel sequence} for $\mathcal{H}$ w.r.t. $\{\mathcal{H}_m\}_{m\in\mathbb{N}}$ and $B$ is called a $g$-Bessel bound. 
  
  If a $g$-frame $\{\Lambda_m\}_{m\in \mathbb{N}}$ ceases to be a $g$-frame when an arbitrary element is removed then $\{\Lambda_m\}_{m\in \mathbb{N}}$ is called a \emph{$g$-exact frame}.
  
  Associated with a $g$-frame $\{\Lambda_m\}_{m\in\mathbb{N}}$, the \emph{$g$-frame operator} $S:\mathcal{H}\rightarrow\mathcal{H}$ is defined by
  \begin{align*}
  	S(h) =\sum_{m\in\mathbb{N}}\Lambda^*_m \Lambda_m h.
  \end{align*}
  The $g$-frame operator $S$ is linear, bounded, self adjoint and invertible.
  
  \begin{defn} \cite{S}
  	Suppose $\{\Lambda_m\}_{m\in\mathbb{N}}$
  	and $\{\Gamma_m\}_{m\in\mathbb{N}}$ are  $g$-frames for $\mathcal{H}$ w.r.t. $\{\mathcal{H}_m\}_{m\in\mathbb{N}}$ such that 
  	\begin{align*}
  		h=\sum_{m\in\mathbb{N}}\Lambda_m^*\Gamma_m h=\sum_{m\in\mathbb{N}}\Gamma_m^*\Lambda_m h,\ \forall h\in\mathcal{H}.
  	\end{align*}
  	Then $\{\Gamma_m\}_{m\in\mathbb{N}}$ is called a \emph{dual $g$-frame} of $\{\Lambda_m\}_{m\in\mathbb{N}}$.
  \end{defn}
  Sun \cite{S} introduced the concept of generalized Riesz basis or $g$-Riesz basis which is the generalization of Riesz basis.
  
  \begin{defn} \cite{S}
  	A sequence $\{\Lambda_m\in\mathcal{L}(\mathcal{H}, \mathcal{H}_m):m \in \mathbb{N}\}$
  	is called a  \emph{generalized Riesz basis} (or \emph{$g$-Riesz basis}) for $\mathcal{H}$ w.r.t. $\{\mathcal{H}_m\}_{m\in\mathbb{N}}$  if 
  	\begin{enumerate}
  		\item [(i)] $ \{\Lambda_m\}_{m\in\mathbb{N}}$ is complete in $\mathcal{H}$ that is $\{h:\Lambda_m h=0, m\in\mathbb{N}\}=\{0\}$
  		\item[(ii)]There exist positive constants $A\leq B$ such that for any finite set $\mathcal{J}\subset\mathbb{N}$ and $h_m\in\mathcal{H}_m$,
  		\begin{align*}
  			A\sum_{m\in\mathcal{J}}\|h_m\|^2\leq\Big\|\sum_{m\in \mathcal{J}}\Lambda_m^*h_m\Big\|^2\leq B\sum_{m\in \mathcal{J}}\|h_m\|^2.
  		\end{align*}
  	\end{enumerate}
 The constants $A$ and $B$ are called lower and upper $g$-Riesz bounds, respectively.
  \end{defn}
  
  \begin{defn} \cite{S}
  	A sequence $\{\Lambda_m\in\mathcal{L}(\mathcal{H}, \mathcal{H}_m):m \in \mathbb{N}\}$
  	is called a \emph{$g$-orthonormal basis} for $\mathcal{H}$ w.r.t. $\{\mathcal{H}_m\}_{m\in\mathbb{N}}$  if 
  	\begin{enumerate}
  		\item[(i)] $\langle \Lambda_{m_1}^*h_{m_1}, \Lambda_{m_2}^*h_{m_2}\rangle=\delta_{m_1,m_2}\langle h_{m_1}, h_{m_2}\rangle $, $\ \forall m_1,m_2\in\mathbb{N},\ h_{m_1}\in\mathcal{H}_{m_1},\ h_{m_2}\in\mathcal{H}_{m_2}$
  		\item[(ii)] 
  		$\sum\limits_{m\in\mathbb{N}}\|\Lambda_mh\|^2=\|h\|^2$ $,\ \forall h\in\mathcal{H}$.
  	\end{enumerate}
  \end{defn} 
  
  Sun \cite{S} characterized $g$-frames, $g$-orthonormal bases and $g$-Riesz bases using orthonormal basis for $\mathcal{H}_m$.
  \begin{thm}\label{th-1s}\cite{S}
  	Let $\Lambda_m\in \mathcal{L}(\mathcal{H}, \mathcal{H}_m)$and $\{e_{n,m}\}_{n\in J_m}$ be an orthonormal basis for $\mathcal{H}_m$, where $J_m\subseteq\mathbb{N}$, $m \in \mathbb{N}$. Then $\{\Lambda_m\}_{m\in\mathbb{N}}$ is a $g$-frame (respectively $g$-Riesz basis, $g$-orthonormal basis) for  $\mathcal{H}$ if and only if $\{\Lambda_m^*e_{n,m}\}_{n\in J_m,m\in\mathbb{N}}$ is a frame (respectively Riesz basis, orthonormal basis) for $\mathcal{H}$.
  \end{thm}

  \section{Weaving generalized frames}
  We begin this section with the definition of weaving $g$-frames in separable Hilbert spaces.
  
  \begin{defn}\cite{KB}
  	Two $g$-frames $\{\Lambda_{m}\}_{m\in\mathbb{N}}$ and $\{\Omega_{m}\}_{m\in\mathbb{N}}$ for $\mathcal{H}$ w.r.t. $\{\mathcal{H}_{m} \}_{m\in\mathbb{N}}$ are called \emph{woven} if there exist positive constants $A\leq B$  such that for any $\sigma\subseteq\mathbb{N}$, $\{\Lambda_{m}\}_{m\in\sigma}\cup\{\Omega_{m}\}_{m\in\sigma^c}$ is a $g$-frame  for $\mathcal{H}$ with lower $g$-frame bound $A$ and upper $g$-frame bound $B$.
  \end{defn}
  The constants $A$ and $B$ are called universal lower $g$-frame bound and universal upper $g$-frame bound, respectively. The supremum of all universal lower $g$-frame bounds is called the \emph{optimal universal lower $g$-frame bound}, and  the infimum of all upper $g$-frame bounds is called the \emph{optimal universal upper $g$-frame bound}.
  
  The following proposition gives the existence of universal  upper $g$-frame bound for any two $g$-frames $\{\Lambda_{m}\}_{m\in\mathbb{N}}$ and $\{\Omega_{m}\}_{m\in\mathbb{N}}$.
  
  \begin{prop}\label{pro1}\cite{VD2}
  	Suppose  $\{\Lambda_{m}\}_{m\in\mathbb{N}}$ and $\{\Omega_{m}\}_{m\in\mathbb{N}}$ are $g$-frames for $\mathcal{H}$ w.r.t. $\{\mathcal{H}_{m}\}_{m\in\mathbb{N}}$ with upper $g$-frame bounds $B_1$ and $B_2$, respectively. Then $B_1+B_2$ is a universal upper $g$-frame bound of $\{\Lambda_{m}\}_{m\in\mathbb{N}}$ and $\{\Omega_{m}\}_{m\in\mathbb{N}}$.
  \end{prop}
  
  In Theorem \ref{th-1s}, $g$-frames are characterized using orthonormal basis for $\mathcal{H}_m$. Following theorem is an extension of Theorem \ref{th-1s} but it characterizes weaving $g$-frames using frames for $\mathcal{H}_m$ instead of orthonormal bases.

  \begin{thm}\label{th-1sd}
  	Let $\{\Lambda_{m}\}_{m\in\mathbb{N}}$ and $\{\Omega_{m}\}_{m\in\mathbb{N}}$ be $g$-frames for $\mathcal{H}$  w.r.t. $\{ \mathcal{H}_m\}_{m\in\mathbb{N}}$. Suppose $\{f_{n,m}\}_{n\in J_m}$ and $\{g_{n,m}\}_{n\in J_m}$ are frames for $\mathcal{H}_m$ with lower frame bounds $A_{1,m}$, $A_{2,m}$ (respectively) and upper frame bounds $B_{1,m}$, $B_{2,m}$ (respectively). If there exist positive constants $A_1<B_1$ and $A_2<B_2$ such that $0<A_1\leq A_{1,m}\leq B_{1,m}\leq B_1<\infty$ and $0<A_2\leq A_{2,m}\leq B_{2,m}\leq B_2<\infty$, for all $m\in\mathbb{N}$, then $\{\Lambda_{m}\}_{m\in\mathbb{N}}$ and $\{\Omega_{m}\}_{m\in\mathbb{N}}$ are weaving $g$-frames for $\mathcal{H}$ if and only if $\{\Lambda_{m}^*f_{n,m}\}_{n\in J_m, m\in\mathbb{N}}$ and $\{\Omega_{m}^*g_{n,m}\}_{n\in J_m, m\in\mathbb{N}}$ are weaving frames for $\mathcal{H}$.
  \end{thm}
  \proof
  First suppose that $\{\Lambda_{m}\}_{m\in\mathbb{N}}$ and $\{\Omega_{m}\}_{m\in\mathbb{N}}$ are weaving $g$-frames for $\mathcal{H}$ with universal lower and upper $g$-frame bounds $A$ and $B$, respectively.

  Let $\sigma$ be any subset of $\mathbb{N}$. Then for any $h\in\mathcal{H}$, we compute
  \begin{align*}
  A\|h\|^2&\leq\sum_{m\in\sigma} \|\Lambda_{m} h\|^2+\sum_{m\in\sigma^c}\|\Omega_m h\|^2\\
  &\leq\sum_{m\in\sigma} \frac{1}{A_{1,m}}\sum_{n\in J_m}|\langle\Lambda_{m} h,f_{n,m}\rangle|^2+\sum_{m\in\sigma^c}\frac{1}{A_{2,m}}\sum_{n\in J_m}|\langle\Omega_m h,g_{n,m}\rangle|^2\\
  &\leq\frac{1}{A_1}\sum_{m\in\sigma}\sum_{n\in J_m}|\langle h,\Lambda_{m}^*f_{n,m}\rangle|^2+\frac{1}{A_{2}}\sum_{m\in\sigma^c}\sum_{n\in J_m}| h,\Omega_m^*g_{n,m}\rangle|^2\\
  &\leq\max\left\{\frac{1}{A_1},\frac{1}{A_{2}}\right\}\left(\sum_{m\in\sigma}\sum_{n\in J_m}|\langle h,\Lambda_{m}^*f_{n,m}\rangle|^2+\sum_{m\in\sigma^c}\sum_{n\in J_m}| h,\Omega_m^*g_{n,m}\rangle|^2\right).
  \end{align*}
  Similarly,
  \begin{align*}
  B\|h\|^2&\geq\sum_{m\in\sigma} \|\Lambda_{m} h\|^2+\sum_{m\in\sigma^c}\|\Omega_m h\|^2\\
  &\geq\sum_{m\in\sigma} \frac{1}{B_{1,m}}\sum_{n\in J_m}|\langle\Lambda_{m} h,f_{n,m}\rangle|^2+\sum_{m\in\sigma^c}\frac{1}{B_{2,m}}\sum_{n\in J_m}|\langle\Omega_m h,g_{n,m}\rangle|^2\\
  &\geq\min\left\{\frac{1}{B_1},\frac{1}{B_{2}}\right\}\left(\sum_{m\in\sigma}\sum_{n\in J_m}|\langle h,\Lambda_{m}^*f_{n,m}\rangle|^2+\sum_{m\in\sigma^c}\sum_{n\in J_m}| h,\Omega_m^*g_{n,m}\rangle|^2\right).
  \end{align*}
  Therefore, $\{\Lambda_{m}^*f_{n,m}\}_{n\in J_m, m\in\mathbb{N}}$ and $\{\Omega_{m}^*g_{n,m}\}_{n\in J_m, m\in\mathbb{N}}$ are weaving frames for $\mathcal{H}$.

  To prove the converse part, suppose $\{\Lambda_{m}^*f_{n,m}\}_{n\in J_m, m\in\mathbb{N}}$ and $\{\Omega_{m}^*g_{n,m}\}_{n\in J_m, m\in\mathbb{N}}$ are weaving frames for $\mathcal{H}$ with universal lower and upper frame bounds $\alpha$ and $\beta$, respectively.
  
  Let $\sigma\subseteq\mathbb{N}$ and $h\in\mathcal{H}$ be arbitrary. Then 
  \begin{align*}
  \alpha\|h\|^2&\leq\sum_{m\in\sigma}\sum_{n\in J_m}|\langle h,\Lambda_{m}^*f_{n,m}\rangle|^2+\sum_{m\in\sigma^c}\sum_{n\in J_m}|\langle h,\Omega_m^*g_{n,m}\rangle|^2\\
  &=\sum_{m\in\sigma}\sum_{n\in J_m}|\langle \Lambda_{m}h,f_{n,m}\rangle|^2+\sum_{m\in\sigma^c}\sum_{n\in J_m}|\langle \Omega_m h,g_{n,m}\rangle|^2\\
  &\leq\sum_{m\in\sigma}B_{1,m}\| \Lambda_{m}h\|^2+\sum_{m\in\sigma^c}B_{2,m}\| \Omega_m h\|^2\\
  &\leq\sum_{m\in\sigma}B_1\| \Lambda_{m}h\|^2+\sum_{m\in\sigma^c}B_2\| \Omega_m h\|^2\\
  &\leq\max\{B_1,B_2\}\left(\sum_{m\in\sigma}\| \Lambda_{m}h\|^2+\sum_{m\in\sigma^c}\| \Omega_m h\|^2\right).
  \end{align*}
  Similarly,
  \begin{align*}
  \beta\|h\|^2&\geq\sum_{m\in\sigma}\sum_{n\in J_m}|\langle h,\Lambda_{m}^*f_{n,m}\rangle|^2+\sum_{m\in\sigma^c}\sum_{n\in J_m}|\langle h,\Omega_m^*g_{n,m}\rangle|^2\\
  &=\sum_{m\in\sigma}\sum_{n\in J_m}|\langle \Lambda_{m}h,f_{n,m}\rangle|^2+\sum_{m\in\sigma^c}\sum_{n\in J_m}|\langle \Omega_m h,g_{n,m}\rangle|^2\\
  &\geq\min\{A_1,A_2\}\left(\sum_{m\in\sigma}\| \Lambda_{m}h\|^2+\sum_{m\in\sigma^c}\| \Omega_m h\|^2\right).
  \end{align*}
  Therefore, $\{\Lambda_{m}\}_{ m\in\mathbb{N}}$ and $\{\Omega_{m}\}_{ m\in\mathbb{N}}$ are weaving $g$-frames for $\mathcal{H}$.
  \endproof
  
  Next example illustrates the above theorem.
  
  \begin{exa}
  Suppose $\mathcal{H}$ is a separable Hilbert space with orthonormal basis $\{e_{n}\}_{n\in\mathbb{N}}$. For $m\in\mathbb{N}$, let $\mathcal{H}_m=\text{span}\{e_m,e_{m+1},e_{m+2}\}$ and define $\Lambda_{m}\in\mathcal{L}(\mathcal{H}, \mathcal{H}_m)$ by
  	\begin{align*}
  \Lambda_{m}(h)&=\langle h,e_m\rangle e_m.
  	\end{align*}
  Here $\Lambda_{m}$ is the orthogonal projection of $\mathcal{H}$ onto $\text{span}\{e_m\}$, so $\Lambda_{m}^*=\Lambda_{m}$, and $\{2e_{m},2e_{m+1},2e_{m+2}\}$ is a frame for $\mathcal{H}_m$ with lower and upper frame bounds both equal to $4$. For any $h\in\mathcal{H}$, we have
  		\begin{align*}
  		\sum_{m\in\mathbb{N}}\|\Lambda_{m}h\|^2=	\sum_{m\in\mathbb{N}}\|\langle h,e_m\rangle e_m\|^2=\sum_{m\in\mathbb{N}}|\langle h,e_m\rangle |^2=\|h\|^2.
  		\end{align*}
  		Therefore, $\{\Lambda_{m}\}_{m\in\mathbb{N}}$ is a $g$-frame for $\mathcal{H}$ w.r.t $\{\mathcal{H}_m\}_{m\in\mathbb{N}}$.\\
  	
  	\vspace{10pt}
  	\textbf{I.}
  	For $m\in\mathbb{N}$, define $\Omega_m\in\mathcal{L}(\mathcal{H}, \mathcal{H}_m)$ by
  	\begin{align*}
  	\Omega_{m}(h)=\begin{cases}\langle h,e_1\rangle e_1+\langle h,e_2\rangle e_2, \quad \text{ if } m=1\\
  	\langle h,e_{m+1}\rangle e_{m+1}, \quad\quad\quad \text{ if } m\geq 2.
  	\end{cases}
  	\end{align*}
  	Since $\Omega_{m}$ is an orthogonal projection, so $\Omega_{m}^*=\Omega_{m}$. Also $\{\Omega_{m}\}_{m\in\mathbb{N}}$ is a $g$-frame for $\mathcal{H}$ w.r.t $\{\mathcal{H}_m\}_{m\in\mathbb{N}}$.\\
  	
  	For $\sigma=\{1\}$ and $h=e_2$, we compute
  		\begin{align*}
  		\sum_{m\in\sigma}(|\langle h,\Lambda_m^*2e_m\rangle|^2&+|\langle h,\Lambda_m^*2e_{m+1}\rangle|^2+|\langle h,\Lambda_m^*2e_{m+2}\rangle|^2)\\
  		&+\sum_{m\in\sigma^c}\left(|\langle h,\Omega_m^*2e_m\rangle|^2+|\langle h,\Omega_m^*2e_{m+1}\rangle|^2+|\langle h,\Omega_m^*2e_{m+2}\rangle|^2\right)\\
  		&=\left(|\langle h,2\Lambda_1e_1\rangle|^2+|\langle h,2\Lambda_1e_{2}\rangle|^2+|\langle h,2\Lambda_1e_{3}\rangle|^2\right)\\
  		&\quad\quad+\sum_{m=2}^\infty\left(|\langle h,2\Omega_me_m\rangle|^2+|\langle h,2\Omega_me_{m+1}\rangle|^2+|\langle h,2\Omega_me_{m+2}\rangle|^2\right)\\	
  		&=|\langle e_2,2e_1\rangle|^2+\sum_{m=2}^\infty|\langle e_2,2e_{m+1}\rangle|^2\\
  		&=0.
  		\end{align*}
  	Thus $\{\Lambda_{m}^*2e_{m},\Lambda_{m}^*2e_{m+1},\Lambda_{m}^*2e_{m+2}\}_{m\in\sigma}\cup\{\Omega_{m}^*2e_{m},\Omega_{m}^*2e_{m+1},\Omega_{m}^*2e_{m+2}\}_{m\in\sigma^c}$ is not a frame for $\mathcal{H}$. Hence $\{\Lambda_{m}^*2e_{m},\Lambda_{m}^*2e_{m+1},\Lambda_{m}^*2e_{m+2}\}_{m\in\mathbb{N}}$ and $\{\Omega_{m}^*2e_{m},\Omega_{m}^*2e_{m+1},\Omega_{m}^*2e_{m+2}\}_{m\in\mathbb{N}}$ are not weaving frames for $\mathcal{H}$, so by Theorem \ref{th-1sd}, $\{\Lambda_{m}\}_{m\in\mathbb{N}}$ and $\{\Omega_{m}\}_{m\in\mathbb{N}}$ are not weaving $g$-frames for $\mathcal{H}$ .\\
  		
  		\vspace{10pt}
  		\textbf{II.}
  			For $m\in\mathbb{N}$, define $\Omega_m\in\mathcal{L}(\mathcal{H}, \mathcal{H}_m)$ by
  		\begin{align*}
  		\Omega_{m}(h)=
  		\langle h,e_{m}\rangle e_{m}+\langle h,e_{m+1}\rangle e_{m+1}.
  		\end{align*}
  Here $\{\Omega_{m}\}_{m\in\mathbb{N}}$ is a $g$-frame for $\mathcal{H}$ w.r.t. $\{\mathcal{H}_m\}_{m\in\mathbb{N}}$ and $\Omega_{m}^*=\Omega_{m}$.\\
  		
  		For $\sigma\subseteq\mathbb{N}$ and $h\in\mathcal{H}$, we have
  		\begin{align*}
  		\sum_{m\in\sigma}(|\langle h,\Lambda_m^*2e_m\rangle|^2&+|\langle h,\Lambda_m^*2e_{m+1}\rangle|^2+|\langle h,\Lambda_m^*2e_{m+2}\rangle|^2)\\
  		&+\sum_{m\in\sigma^c}\left(|\langle h,\Omega_m^*2e_m\rangle|^2+|\langle h,\Omega_m^*2e_{m+1}\rangle|^2+|\langle h,\Omega_m^*2e_{m+2}\rangle|^2\right)\\
  			&=\sum_{m\in\sigma}|\langle h,2e_m\rangle|^2
  		+\sum_{m\in\sigma^c}(|\langle h,2e_m\rangle|^2+|\langle h,2e_{m+1}\rangle|^2)\\
  		&=4\|h\|^2+4\sum_{m\in\sigma^c}|\langle h,e_{m+1}\rangle|^2
  		\end{align*}	
  	Therefore	
 \begin{align*}
  	 \|h\|^2&\leq\sum_{m\in\sigma}(|\langle h,\Lambda_m^*2e_m\rangle|^2+|\langle h,\Lambda_m^*2e_{m+1}\rangle|^2+|\langle h,\Lambda_m^*2e_{m+2}\rangle|^2)\\
  		+\sum_{m\in\sigma^c}&\left(|\langle h,\Omega_m^*2e_m\rangle|^2+|\langle h,\Omega_m^*2e_{m+1}\rangle|^2+|\langle h,\Omega_m^*2e_{m+2}\rangle|^2\right)\leq 8\|h\|^2.
  		\end{align*}
  		Hence, $\{\Lambda_{m}^*2e_{m},\Lambda_{m}^*2e_{m+1},\Lambda_{m}^*2e_{m+2}\}_{m\in\mathbb{N}}$ and $\{\Omega_{m}^*2e_{m},\Omega_{m}^*2e_{m+1},\Omega_{m}^*2e_{m+2}\}_{m\in\mathbb{N}}$ are weaving frames for $\mathcal{H}$, so by Theorem \ref{th-1sd}, $\{\Lambda_{m}\}_{m\in\mathbb{N}}$ and $\{\Omega_{m}\}_{m\in\mathbb{N}}$ are weaving $g$-frames for $\mathcal{H}$ .
   \end{exa}
 In the following theorem, we show the relation between optimal $g$-frame bounds and optimal universal $g$-frame bounds.
 
 \begin{thm}\label{th-3.2}
 	Let   $\{\Lambda_{m}\}_{m\in\mathbb{N}}$ and $\{\Omega_{m}\}_{m\in\mathbb{N}}$ be $g$-frames  for $\mathcal{H}$ w.r.t. $\{ \mathcal{H}_m\}_{m\in\mathbb{N}}$ with optimal lower $g$-frame bounds $A_1$, $A_2$ (respectively) and optimal upper $g$-frame bounds $B_1$, $B_2$ (respectively). If $\{\Lambda_{m}\}_{m\in\mathbb{N}}$ and $\{\Omega_{m}\}_{m\in\mathbb{N}}$ are weaving $g$-frames for $\mathcal{H}$  with optimal universal lower and upper  $g$-frame bounds $A$ and $B$, respectively, then $A\leq\min\{A_1,A_2\}$ and $B\geq\max\{B_1,B_2\}$.
 \end{thm}
 \proof
Choose $\sigma=\mathbb{N}$. Then for any $h\in\mathcal{H}$, we have
 \begin{align*}
	A\|h\|^2\leq\sum_{m\in\sigma} \|\Lambda_{m} h\|^2+\sum_{m\in\sigma^c}\|\Omega_m h\|^2=\sum_{m\in\mathbb{N}} \|\Lambda_{m} h\|^2\leq B\|h\|^2
\end{align*}
Thus $A$ and $B$ are lower and upper $g$-frame bounds, respectively of  $\{\Lambda_{m}\}_{m\in\mathbb{N}}$. Since  $A_1$ and $B_1$ are optimal $g$-frame bounds of  $\{\Lambda_{m}\}_{m\in\mathbb{N}}$, $A\leq A_1$ and $B\geq B_1$.

Similarly  $A\leq A_2$ and $B\geq B_2$. Hence $A\leq\min\{A_1,A_2\}$ and $B\geq\max\{B_1,B_2\}$.
 \endproof
 
 Let us illustrate an example of Theorem \ref{th-3.2} where strict inequalities follows.
 
 \begin{exa}
 	Suppose $\mathcal{H}$ is a separable Hilbert space with orthonormal basis $\{e_{n}\}_{n\in\mathbb{N}}$. Let $\mathcal{H}_1=\text{span}\{e_1\}$, $\mathcal{H}_2=\mathcal{H}_3=\text{span}\{e_2\}$, $\mathcal{H}_4=\mathcal{H}_5=\text{span}\{e_3\}$ and for $m\in\mathbb{N}\setminus\{1,2,3,4,5\}$, let $\mathcal{H}_m=\text{span}\{e_{m-2}\}$.
 	
For $m\in\mathbb{N}\setminus\{2,3,4,5\}$, $\Lambda_{m}$ and $\Omega_{m}$ are orthogonal projection of $\mathcal{H}$ onto $\mathcal{H}_m$. Define $\Lambda_{2},\Lambda_{3},\Lambda_{4}, \Lambda_{5}, \Omega_{2},  \Omega_{3}, \Omega_{4}, \Omega_{5}$ as
 	\begin{align*}
 		\Lambda_{2}(h)&=\langle h,\frac{e_2}{\sqrt{2}}\rangle e_2,\ \Lambda_{3}(h)=\langle h,\frac{e_2}{\sqrt{2}}\rangle e_2,\ \Lambda_{4}(h)=\langle h,e_3\rangle e_3,\ \Lambda_{5}(h)=0\\
 		\Omega_{2}(h)&=\langle h,e_2\rangle e_2,\ \Omega_{3}(h)=0,\ \Omega_{4}(h)=\langle h,\frac{e_3}{\sqrt{2}}\rangle e_3,\ \Omega_{5}(h)=\langle h,\frac{e_3}{\sqrt{2}}\rangle e_3
 	\end{align*}
 	 For any $h\in\mathcal{H}$, we have
 	\begin{align*}
 		\sum_{m\in\mathbb{N}}\|\Lambda_{m}h\|^2=\sum_{m\in\mathbb{N}}\|\langle h,e_m\rangle\|^2=	\sum_{m\in\mathbb{N}}\|\Omega_{m}h\|^2=\|h\|^2.
 	\end{align*}
 	Therefore, $\{\Lambda_{m}\}_{m\in\mathbb{N}}$ and $\{\Omega_{m}\}_{m\in\mathbb{N}}$ are $g$-frames for $\mathcal{H}$ w.r.t $\{\mathcal{H}_m\}_{m\in\mathbb{N}}$ with optimal lower and upper $g$-frame bounds equal to $1$. Thus $A_1=A_2=B_1=B_2=1$.\\
 	For any $\sigma\subset\mathbb{N}$ and for any $h\in\mathcal{H}$,
 	 	\begin{align*}
 	 	\frac{1}{2}\|h\|^2&=\frac{1}{2}\sum_{m\in\mathbb{N}}\|\langle h,e_m\rangle\|^2\\
 	 	&\leq\sum_{m\in\sigma}\|\Lambda_{m}h\|^2+\sum_{m\in\sigma^c}\|\Omega_{m}h\|^2\\
 	 	&\leq\frac{3}{2}\sum_{m\in\mathbb{N}}\|\langle h,e_m\rangle\|^2\\
 	 	&=\frac{3}{2}\|h\|^2.
 	 \end{align*}
 Thus,  $\{\Lambda_{m}\}_{m\in\mathbb{N}}$ and $\{\Omega_{m}\}_{m\in\mathbb{N}}$ are weaving $g$-frames for  $\mathcal{H}$ with universal lower and upper $g$-frame bounds $\frac{1}{2}$ and $\frac{3}{2}$, respectively.
 
 For $\sigma_1=\{2\}$, $\sigma_2=\{4\}$, $h_1=e_2$, $h_2=e_3$, we have
 \begin{align*}
 	\frac{1}{2}\|h_1\|^2&=\sum_{m\in\sigma_1}\|\Lambda_{m}h\|^2+\sum_{m\in\sigma_1^c}\|\Omega_{m}h\|^2\\
	\frac{3}{2}\|h_2\|^2&=\sum_{m\in\sigma_2}\|\Lambda_{m}h_2\|^2+\sum_{m\in\sigma_2^c}\|\Omega_{m}h_2\|^2.
\end{align*}
Therefore, optimal universal upper $g$-frame bound is $\frac{3}{2}$ and optimal universal lower $g$-frame bound is $\frac{1}{2}$. Hence $A=\frac{1}{2}<\min\{A_1,A_2\}$ and $B=\frac{3}{2}>\max\{B_1,B_2\}$.
 \end{exa}
  
  In the next theorem, we show that the sum of the optimal $g$-frame bounds of two weaving $g$-frames are never the optimal universal $g$-frame bounds.

  \begin{thm}
  	Let   $\{\Lambda_{m}\}_{m\in\mathbb{N}}$ and $\{\Omega_{m}\}_{m\in\mathbb{N}}$ be $g$-frames  for $\mathcal{H}$ w.r.t. $\{ \mathcal{H}_m\}_{m\in\mathbb{N}}$ with optimal lower $g$-frame bounds $A_1$, $A_2$ (respectively) and optimal upper $g$-frame bounds $B_1$, $B_2$ (respectively). If $\{\Lambda_{m}\}_{m\in\mathbb{N}}$ and $\{\Omega_{m}\}_{m\in\mathbb{N}}$ are weaving $g$-frames for $\mathcal{H}$, then $A_1+A_2$ is not the optimal universal lower $g$-frame bound and $B_1+B_2$ is not the optimal universal upper $g$-frame bound.
  \end{thm}
  \proof
  Suppose $\{\Lambda_{m}\}_{m\in\mathbb{N}}$ and $\{\Omega_{m}\}_{m\in\mathbb{N}}$ are weaving $g$-frames for $\mathcal{H}$ with universal lower and upper $g$-frame bounds $A$ and $B$, respectively. Since $\min\{A_1,A_2\}<A_1+A_2$, by Theorem \ref{th-3.2} optimal universal lower $g$-frame bound is not equal to $A_1+A_2$. 
  
  Let if possible suppose that $B_1+B_2$ is the optimal universal upper $g$-frame bound of weaving $g$-frames $\{\Lambda_{m}\}_{m\in\mathbb{N}}$ and $\{\Omega_{m}\}_{m\in\mathbb{N}}$. For any $\epsilon>0$, $B_1+B_2-\epsilon$ is not a universal upper $g$-frame bound. Thus there exists $\sigma\subseteq\mathbb{N}$ and $h\in \mathcal{H}$ such that 
  \begin{align*}
  \sum\limits_{m\in\sigma} \|\Lambda_{m} h\|^2+\sum\limits_{m\in\sigma^c}\|\Omega_m h\|^2>(B_1+B_2-\epsilon)\|h\|^2.
  \end{align*}
  Take $h_1=\frac{h}{\|h\|}$. Then
  \begin{align*}
  \sum\limits_{m\in\sigma} \|\Lambda_{m} h_1\|^2+\sum\limits_{m\in\sigma^c}\|\Omega_m h_1\|^2>B_1+B_2-\epsilon\geq \sum\limits_{m\in\mathbb{N}} \|\Lambda_{m} h_1\|^2+\sum\limits_{m\in\mathbb{N}}\|\Omega_m h_1\|^2-\epsilon
\end{align*}
Thus, we have
  \begin{align*}
  A\leq\sum\limits_{m\in\sigma^c} \|\Lambda_{m} h_1\|^2+\sum\limits_{m\in\sigma}\|\Omega_m h_1\|^2<\epsilon
  \end{align*}
  Since $\epsilon>0$ was arbitrary, so $A=0$, a contradiction. Therefore, $B_1+B_2$ is not the optimal universal upper $g$-frame bound.
  \endproof
  Dual $g$-frames provide the series representation of each vector in $\mathcal{H}$ in terms of $g$-frame elements. Following theorem shows that a $g$-frame and its dual $g$-frame are woven.
  
  \begin{thm}\label{th*}
  	Let   $\{\Lambda_{m}\}_{m\in\mathbb{N}}$ be a $g$-frame for $\mathcal{H}$ w.r.t. $\{ \mathcal{H}_m\}_{m\in\mathbb{N}}$ with upper $g$-frame bound $B_1$ and let $\{\Gamma_{m}\}_{m\in\mathbb{N}}$ be its dual $g$-frame with upper $g$-frame bound $B_2$. Then $\{\Lambda_{m}\}_{m\in\mathbb{N}}$ and $\{\Gamma_{m}\}_{m\in\mathbb{N}}$ are weaving $g$-frames for $\mathcal{H}$ with universal lower $g$-frame bound $\min\left\{\frac{1}{2B_1},\frac{1}{2B_2}\right\}$ and universal upper $g$-frame bound $B_1+B_2$.
  \end{thm}
  \proof
  Let $\sigma$ be any subset of $\mathbb{N}$. By using Cauchy Schwartz inequality, we compute
  \begin{align*}
  \|h\|^4&=|\langle h,h\rangle|^2\\
  &=|\langle \sum_{m\in\mathbb{N}}\Lambda_{m}^* \Gamma_m h,h\rangle|^2\\
  &=|\sum_{m\in\sigma}\langle \Lambda_{m}^* \Gamma_m h,h\rangle+\sum_{m\in\sigma^c}\langle \Lambda_{m}^* \Gamma_m h,h\rangle|^2\\
  &=|\sum_{m\in\sigma}\langle \Gamma_m h,\Lambda_{m}h\rangle+\sum_{m\in\sigma^c}\langle  \Gamma_m h,\Lambda_{m}h\rangle|^2\\
  &\leq2|\sum_{m\in\sigma}\langle \Gamma_m h,\Lambda_{m}h\rangle|^2+2|\sum_{m\in\sigma^c}\langle  \Gamma_m h,\Lambda_{m}h\rangle|^2\\
 \quad\quad\quad   &\leq2\sum_{m\in\sigma}\|\Gamma_m h\|^2\sum_{m\in\sigma} \|\Lambda_{m}h\|^2+2\sum_{m\in\sigma^c}\|\Gamma_m h\|^2\sum_{m\in\sigma^c} \|\Lambda_{m}h\|^2\\
  &\leq2B_2\|h\|^2\sum_{m\in\sigma} \|\Lambda_{m}h\|^2+2B_1\|h\|^2\sum_{m\in\sigma^c}\|\Gamma_m h\|^2\\
  &\leq\max\{2B_1,2B_2\}\|h\|^2\left(\sum_{m\in\sigma} \|\Lambda_{m}h\|^2+\sum_{m\in\sigma^c}\|\Gamma_m h\|^2\right), \ \forall h\in\mathcal{H}.
\end{align*}
Therefore
\begin{align*}
  \min\left\{\frac{1}{2B_1},\frac{1}{2B_2}\right\}&\|h\|^2\leq \sum\limits_{m\in\sigma} \|\Lambda_{m}h\|^2+\sum\limits_{m\in\sigma^c}\|\Gamma_m h\|^2\leq (B_1+B_2)\|h\|^2,\ \forall h\in\mathcal{H}.
  \end{align*}
  \endproof
  
  Since the $g$-frame operator and its inverse are self adjoint and positive, so their square roots exist. In the next theorem, we construct a new family of weaving $g$-frames using the existing family of weaving $g$-frames and the square root of $S^{-1}$. 
  
  \begin{thm}
  	Let $\{\Lambda_{m}\}_{m\in\mathbb{N}}$ and $\{\Omega_{m}\}_{m\in\mathbb{N}}$ be weaving $g$-frames  for $\mathcal{H}$ w.r.t. $\mathcal{H}_m$. If $S$ is the $g$-frame operator of $\{\Lambda_{m}\}_{m\in\mathbb{N}}$, then $\{\Lambda_{m}S^{-\frac{1}{2}}\}_{m\in\mathbb{N}}$ and $\{\Omega_{m}S^{-\frac{1}{2}}\}_{m\in\mathbb{N}}$ are weaving $g$-frames for $\mathcal{H}$  w.r.t. $\mathcal{H}_m$.
  \end{thm}
  \proof
   Suppose  $\{\Lambda_{m}\}_{m\in\mathbb{N}}$ and $\{\Omega_{m}\}_{m\in\mathbb{N}}$ are weaving $g$-frames  for $\mathcal{H}$ with universal lower $g$-frame bound $A$ and universal upper $g$-frame bound $B$. Then $\{\Lambda_{m}\}_{m\in\mathbb{N}}$ is a $g$-frame  for $\mathcal{H}$ with lower $g$-frame bound $A$ and upper $g$-frame bound $B$, and hence $B^{-1}I\leq S^{-1}\leq A^{-1}I$, where $I$ is the identity operator on $\mathcal{H}$.
  
  For any subset $\sigma$ of $\mathbb{N}$ and $h\in\mathcal{H}$, we compute
  \begin{align*}
  	\sum\limits_{m\in\sigma} \|\Lambda_{m}S^{-\frac{1}{2}}h\|^2+\sum\limits_{m\in\sigma^c}\|\Omega_mS^{-\frac{1}{2}} h\|^2&\geq A\|S^{-\frac{1}{2}}h\|^2\\
  	&= A\langle S^{-\frac{1}{2}}h,S^{-\frac{1}{2}}h\rangle\\
  	&= A\langle (S^{-\frac{1}{2}})^*S^{-\frac{1}{2}}h,h\rangle\\
  	&= A\langle S^{-\frac{1}{2}}S^{-\frac{1}{2}}h,h\rangle\\
  	&= A\langle S^{-1}h,h\rangle\\
  	&\geq\frac{A}{B}\|h\|^2.
  \end{align*}
  
  For the universal upper $g$-frame bound, we have
  \begin{align*}
  	\sum\limits_{m\in\sigma} \|\Lambda_{m}S^{-\frac{1}{2}}h\|^2+\sum\limits_{m\in\sigma^c}\|\Omega_mS^{-\frac{1}{2}} h\|^2&\leq B\|S^{-\frac{1}{2}}h\|^2\\
  	&= B\langle S^{-\frac{1}{2}}h,S^{-\frac{1}{2}}h\rangle\\
  	&= B\langle S^{-1}h,h\rangle\\
  	&\leq\frac{B}{A}\|h\|^2.
  \end{align*}
  Therefore, $\{\Lambda_{m}S^{-\frac{1}{2}}\}_{m\in\mathbb{N}}$ and $\{\Omega_{m}S^{-\frac{1}{2}}\}_{m\in\mathbb{N}}$ are weaving $g$-frames for $\mathcal{H}$ with universal lower $g$-frame bound $\frac{A}{B}$ and universal upper $g$-frame bound $\frac{B}{A}$.
  \endproof

  \section{Weaving generalized Riesz bases}
  We start this section with the definition of weaving $g$-Riesz bases and weaving $g$-orthonormal bases.
  
  \begin{defn}
  	Two $g$-Riesz bases $\{\Lambda_{m}\}_{m\in\mathbb{N}}$ and $\{\Omega_{m}\}_{m\in\mathbb{N}}$ for $\mathcal{H}$ w.r.t. $\{\mathcal{H}_{m} \}_{m\in\mathbb{N}}$ are called \emph{weaving $g$-Riesz bases} if there exist positive constants $A\leq B$  such that for any $\sigma\subseteq\mathbb{N}$, $\{\Lambda_{m}\}_{m\in\sigma}\cup\{\Omega_{m}\}_{m\in\sigma^c}$ is a $g$-Riesz basis  for $\mathcal{H}$ with lower $g$-Riesz bound $A$ and upper $g$-Riesz bound $B$.
  \end{defn}
  The constants $A$ and $B$ are called universal lower $g$-Riesz bound and universal upper $g$-Riesz bound, respectively. 
  
  \begin{defn}
  	Two $g$-orthonormal bases $\{\Lambda_{m}\}_{m\in\mathbb{N}}$ and $\{\Omega_{m}\}_{m\in\mathbb{N}}$ for $\mathcal{H}$ w.r.t. $\{\mathcal{H}_{m} \}_{m\in\mathbb{N}}$ are called \emph{weaving $g$-orthonormal bases} if for any $\sigma\subseteq\mathbb{N}$, $\{\Lambda_{m}\}_{m\in\sigma}\cup\{\Omega_{m}\}_{m\in\sigma^c}$ is a $g$-orthonormal basis  for $\mathcal{H}$.
  \end{defn}

  In the next theorem, we show that weaving $g$-orthonormal bases remains to be woven even after applying unitary operator.
  
  \begin{thm}\label{th-1d}
  	If   $\{\Lambda_{m}\}_{m\in\mathbb{N}}$ and $\{\Omega_{m}\}_{m\in\mathbb{N}}$ are weaving $g$-orthonormal bases for $\mathcal{H}$ w.r.t. $\{\mathcal{H}_{m}\}_{m\in\mathbb{N}}$ and $U$ is any unitary operator on $\mathcal{H}$, then $\{\Lambda_{m}U\}_{m\in\mathbb{N}}$ and $\{\Omega_{m}U\}_{m\in\mathbb{N}}$ are weaving $g$-orthonormal bases for $\mathcal{H}$.
  \end{thm}
  \proof
  Suppose $\sigma$ is any subset of $\mathbb{N}$. Then
  \begin{align*}
  \sum\limits_{m\in\sigma} \|\Lambda_{m} Uh\|^2+\sum\limits_{m\in\sigma^c}\|\Omega_mU h\|^2=\|Uh\|^2=\|h\|^2,\ \forall h\in\mathcal{H}.
  \end{align*}
   Suppose $m_1, m_2\in\mathbb{N}$ are arbitrary. Then for any $h_{m_1}\in\mathcal{H}_{m_1}$ and $h_{m_2}\in\mathcal{H}_{m_2}$, we have
  \begin{align*}
  \langle(\Lambda_{m_1}U)^*h_{m_1},(\Omega_{m_2}U)^*h_{m_2}\rangle=\langle U^*\Lambda_{m_1}^*h_{m_1},U^*\Omega_{m_2}^*h_{m_2}\rangle=\langle \Lambda_{m_1}^*h_{m_1},\Omega_{m_2}^*h_{m_2}\rangle=\delta_{m_1,m_2}\langle h_{m_1},h_{m_2}\rangle\\
  \text{if } m_1\in\sigma,\ m_2\in\sigma^c,\\
  \langle(\Lambda_{m_1}U)^*h_{m_1},(\Lambda_{m_2}U)^*h_{m_2}\rangle=\langle U^*\Lambda_{m_1}^*h_{m_1},U^*\Lambda_{m_2}^*h_{m_2}\rangle=\langle \Lambda_{m_1}^*h_{m_1},\Lambda_{m_2}^*h_{m_2}\rangle=\delta_{m_1,m_2}\langle h_{m_1},h_{m_2}\rangle\\
  \text{if } m_1,\ m_2 \in\sigma,\\
  \langle(\Omega_{m_1}U)^*h_{m_1},(\Omega_{m_2}U)^*h_{m_2}\rangle=\langle U^*\Omega_{m_1}^*h_{m_1},U^*\Omega_{m_2}^*h_{m_2}\rangle=\langle \Omega_{m_1}^*h_{m_1},\Omega_{m_2}^*h_{m_2}\rangle=\delta_{m_1,m_2}\langle h_{m_1},h_{m_2}\rangle\\
  \text{if } m_1, \ m_2\in\sigma^c.
  \end{align*}
  Therefore, $\{\Lambda_{m}U\}_{m\in\mathbb{N}}$ and $\{\Omega_{m}U\}_{m\in\mathbb{N}}$ are weaving $g$-orthonormal bases for $\mathcal{H}$.
  \endproof
  
  We obtain the following corollary to the above theorem for $g$-orthonormal basis.
  \begin{cor}\label{cr-1d}
  	If  $\{\Lambda_{m}\}_{m\in\mathbb{N}}$ is a $g$-orthonormal basis for $\mathcal{H}$ w.r.t. $\{\mathcal{H}_{m}\}_{m\in\mathbb{N}}$ and $U$ is any unitary operator on $\mathcal{H}$, then $\{\Lambda_{m}U\}_{m\in\mathbb{N}}$ is a $g$-orthonormal basis for $\mathcal{H}$.
  \end{cor}
  
  Unitary operators are surjective isometries. Both surjectivity and isomerty are necessary in Theorem \ref{th-1d}  and Corollary \ref{cr-1d} is justified in the following example.
  \begin{exa}
  Suppose $\mathcal{H}=\ell^2(\mathbb{N})$ with canonical orthonormal basis $\{e_{n}\}_{n\in\mathbb{N}}$. For $m\in\mathbb{N}$, let $\mathcal{H}_m=\text{span}\{e_m\}$ and define $\Lambda_{m}\in\mathcal{L}(\mathcal{H}, \mathcal{H}_m)$ by
  	\begin{align*}
  \Lambda_{m}(h)=\langle h,e_m\rangle e_m.
  	\end{align*}
  	Here $\Lambda_{m}$ is the orthogonal projection of $\mathcal{H}$ onto $\mathcal{H}_m$, so $\Lambda_{m}^*=\Lambda_{m}$.
  	
  	For any $h\in\mathcal{H}$, we have
  	\begin{align*}
  	\sum_{m\in\mathbb{N}}\|\Lambda_{m}h\|^2=	\sum_{m\in\mathbb{N}}\|\langle h,e_m\rangle e_m\|^2=\sum_{m\in\mathbb{N}}|\langle h,e_m\rangle |^2=\|h\|^2.
  	\end{align*}
  Suppose $m_1, m_2\in\mathbb{N}$ are arbitrary. Then for any $h_{m_1}\in\mathcal{H}_{m_1}$ and $h_{m_2}\in\mathcal{H}_{m_2}$, we have
  	\begin{align*}
  	\langle\Lambda_{m_1}^*h_{m_1},\Lambda_{m_2}^*h_{m_2}\rangle&=\langle\Lambda_{m_1}h_{m_1},\Lambda_{m_2}h_{m_2}\rangle\\
  	&=\langle\langle h_{m_1},e_{m_1}\rangle e_{m_1},\langle h_{m_2},e_{m_2}\rangle e_{m_2}\rangle\\
  	&=\langle h_{m_1},e_{m_1}\rangle\langle e_{m
  	_2}, h_{m_2}\rangle\langle e_{m_1}, e_{m_2}\rangle\\
  	&=\langle h_{m_1},e_{m_1}\rangle\langle e_{m
  		_2}, h_{m_2}\rangle\langle e_{m_1}, e_{m_2}\rangle\langle e_{m_1}, e_{m_2}\rangle\\
  	&=\delta_{m_1,m_2}\langle h_{m_1},e_{m_1}\rangle\langle e_{m_2}, h_{m_2}\rangle\langle e_{m_1}, e_{m_2}\rangle,\\
  	\langle h_{m_1}, h_{m_2}\rangle&=\langle\langle h_{m_1},e_{m_1}\rangle e_{m_1},\langle h_{m_2},e_{m_2}\rangle e_{m_2}\rangle\\
  	&=\langle h_{m_1},e_{m_1}\rangle\langle e_{m_2}, h_{m_2}\rangle\langle e_{m_1}, e_{m_2}\rangle.
  	\end{align*}
  	Thus $	\langle\Lambda_{m_1}^*h_{m_1},\Lambda_{m_2}^*h_{m_2}\rangle=\delta_{m_1,m_2}\langle h_{m_1}, h_{m_2}\rangle$ and hence, $\{\Lambda_{m}\}_{m\in\mathbb{N}}$ is a $g$-orthonormal basis for $\mathcal{H}$ w.r.t. $\{\mathcal{H}_m\}_{m\in\mathbb{N}}$.\\
  	
  	\vspace{10pt}
  	\textbf{I.}
  	Define $U:\mathcal{H}\rightarrow\mathcal{H}$ by $U(h)=2h$. Then $U$ is a bounded, linear and surjective operator but $U$ is not an isometry. For $h\in\mathcal{H}$, we have
  		\begin{align*}
  		\sum_{m\in\mathbb{N}}\|\Lambda_{m}Uh\|^2=\sum_{m\in\mathbb{N}}\|\Lambda_{m}2h\|^2=	\sum_{m\in\mathbb{N}}\|\langle 2h,e_m\rangle e_m\|^2=\sum_{m\in\mathbb{N}}|\langle 2h,e_m\rangle |^2=4\|h\|^2.
  		\end{align*}
  	Therefore, $\{\Lambda_{m}U\}_{m\in\mathbb{N}}$ is not a $g$-orthonormal basis for $\mathcal{H}$. Hence isometry of $U$ is necessary in Corollary \ref{cr-1d}.
  	
  	If $\Omega_{m}=\Lambda_m$, $m\in\mathbb{N}$,  then $\{\Lambda_{m}\}_{m\in\mathbb{N}}$ and $\{\Omega_{m}\}_{m\in\mathbb{N}}$ are weaving $g$-orthonormal bases for $\mathcal{H}$. But $\{\Lambda_{m}U\}_{m\in\mathbb{N}}$ and $\{\Omega_{m}U\}_{m\in\mathbb{N}}$ are not weaving $g$-orthonormal bases for $\mathcal{H}$ as $\{\Lambda_{m}U\}_{m\in\mathbb{N}}$ is not a $g$-orthonormal basis for $\mathcal{H}$. Hence isometry of $U$ is necessary in Theorem \ref{th-1d}.\\
  		
  		\vspace{10pt}
  		\textbf{II.}
  		Define $U:\mathcal{H}\rightarrow\mathcal{H}$ by $U(a_1,a_2,a_3,a_4,\ldots)=(0,a_1,a_2,a_3,a_4,\ldots)$. Then $U$ is a bounded, linear and isometry operator but $U$ is not surjective. Here $U$ is a right shift operator, so $U^*$ is the left shift operator.
  		
  		Since $\{e_m\}$ is an orthonormal basis of $\mathcal{H}_m$ and $\|(\Lambda_1U)^*e_1\|=\|U^*\Lambda_1e_1\|=\|U^*e_1\|=\|(0,0,0,\ldots)\|=0$, so $\{(\Lambda_m U)^*e_m\}_{m\in\mathbb{N}}$ is not an orthonormal basis for $\mathcal{H}$. By Theorem \ref{th-1s}, $\{\Lambda_{m}U\}_{m\in\mathbb{N}}$ is not a $g$-orthonormal basis for $\mathcal{H}$. Hence surjectivity of $U$ is necessary in Corollary \ref{cr-1d}.
  			
  			If $\Omega_{m}=\Lambda_m$, $m\in\mathbb{N}$,  then $\{\Lambda_{m}\}_{m\in\mathbb{N}}$ and $\{\Omega_{m}\}_{m\in\mathbb{N}}$ are weaving $g$-orthonormal bases for $\mathcal{H}$. But $\{\Lambda_{m}U\}_{m\in\mathbb{N}}$ and $\{\Omega_{m}U\}_{m\in\mathbb{N}}$ are not weaving $g$-orthonormal bases for $\mathcal{H}$. Hence surjectivity of $U$ is necessary in Theorem \ref{th-1d}.\\
   \end{exa}
  $G$-Riesz bases are generalization of Riesz bases but still some properties of weaving $g$-Riesz bases and weaving Riesz bases are different. Next two examples highlight these differences.
  
 It is shown in Theorem \ref{th-1c} that if a frame and a Riesz basis are woven, then the frame must be a Riesz basis. However this is not in the case of $g$-Riesz basis. 
  
  \begin{remark}If a $g$-frame and a $g$-Riesz basis are woven then the $g$-frame need not be a $g$-Riesz basis. 
  \end{remark}
  Next example justifies the above remark.
  \begin{exa}\label{eg4.4}
  	Suppose $\mathcal{H}$ is a separable Hilbert space with orthonormal basis $\{e_{n,m}\}_{n,m\in\mathbb{N}}$. For $m\in\mathbb{N}$, define $\Lambda_{m},\Omega_{m}:\mathcal{H}\rightarrow\ell^2(\mathbb{N})$ by
  	\begin{align*}
  &\Lambda_{m}(h)=\begin{cases}
  \{\langle h,e_{n,k}\rangle\}_{k\in\mathbb{N}}\quad \text{ if }m=2n\\
  \{\langle h,e_{n,k}\rangle\}_{k\in\mathbb{N}}\quad \text{ if }m=2n-1
  \end{cases}\\
  &\Omega_{m}(h)=\begin{cases}
  \{\langle h,e_{n,2k}\rangle\}_{k\in\mathbb{N}}\quad \quad \text{ if }m=2n\\
  \{\langle h,e_{n,2k-1}\rangle\}_{k\in\mathbb{N}}\quad \text{ if }m=2n-1.
  \end{cases}
  	\end{align*}
  	For any $h\in\mathcal{H}$, we have
  	\begin{align*}
  	\sum_{m\in\mathbb{N}}\|\Lambda_{m}h\|^2=2\sum_{n\in\mathbb{N}}\|\{\langle h,e_{n,k}\rangle\}_{k\in\mathbb{N}}\|^2=2\sum_{n\in\mathbb{N}}\sum_{k\in\mathbb{N}}|\langle h,e_{n,k}\rangle|^2=2\|h\|^2.
  	\end{align*}
  	Thus $\{\Lambda_{m}\}_{m\in\mathbb{N}}$ is a $g$-frame for $\mathcal{H}$. Since $\{\Lambda_{m}\}_{m\in\mathbb{N}\setminus\{2\}}$ is also a $g$-frame for $\mathcal{H}$, so $\{\Lambda_{m}\}_{m\in\mathbb{N}}$ is not a $g$-exact frame and hence it is not a $g$-Riesz basis for $\mathcal{H}$ as every $g$-Riesz basis is a $g$-exact frame.
  	
  	Let $\{e_i\}_{i\in\mathbb{N}}$ be the canonical orthonormal basis for $\ell^2(\mathbb{N})$. Then
  	\begin{align*}
  	\Omega_{m}^*(e_i)=\begin{cases}
  	e_{n,2i-1}\quad &\text{if } m=2n-1\\
  	e_{n,2i} &\text{if } m=2n.
  	\end{cases}
  	\end{align*}
  	Since $\{\Omega_{m}^*e_{i}\}_{i\in\mathbb{N},m\in\mathbb{N}}=\{e_{i,m}\}_{i\in\mathbb{N},m\in\mathbb{N}}$ is a Riesz basis for $\mathcal{H}$ being an orthonormal basis, so by Theorem \ref{th-1s}, $\{\Omega_{m}\}_{m\in\mathbb{N}}$ is a $g$-Riesz basis for $\mathcal{H}$.
  	
  	For any $\sigma\subseteq\mathbb{N}$ and $h\in\mathcal{H}$, we have
  	\begin{align*}
  	\|h\|^2&=\sum_{m\in\mathbb{N}}\sum_{n\in\mathbb{N}}|\langle h,e_{n,m}\rangle|^2\\
  	&\leq\sum_{m\in\sigma}\|\Lambda_{m}h\|^2+\sum_{m\in\sigma^c}\|\Omega_{m}h\|^2\\
  	&\leq2\sum_{m\in\mathbb{N}}\sum_{n\in\mathbb{N}}|\langle h,e_{n,m}\rangle|^2=2\|h\|^2.
  	\end{align*}
Therefore, $\{\Lambda_{m}\}_{m\in\mathbb{N}}$ and $\{\Omega_{m}\}_{m\in\mathbb{N}}$ are weaving $g$-frames for $\mathcal{H}$, where $\{\Omega_{m}\}_{m\in\mathbb{N}}$ is a $g$-Riesz basis and $\{\Lambda_{m}\}_{m\in\mathbb{N}}$ is not a $g$-Riesz basis.
  \end{exa}
 It is presented in Theorem \ref{th-2c} that if two Riesz bases are woven, then every weaving is a Riesz basis. Since exact frames are same as Riesz bases, so conclusion of Theorem \ref{th-2c} also holds for exact frames i.e. if two exact frames are woven, then every weaving is a Riesz basis (or an exact frame). But this is not true for  $g$-exact frames.
  
  \begin{remark} If two $g$-exact frames $\{\Lambda_{m}\}_{m\in\mathbb{N}}$ and $\{\Omega_{m}\}_{m\in\mathbb{N}}$ are weaving $g$-frames for $\mathcal{H}$, then for $\sigma\subset\mathbb{N}$, $\{\Lambda_{m}\}_{m\in\sigma}\cup\{\Omega_{m}\}_{m\in\sigma^c}$ need not be a $g$-exact frame and hence need not be a $g$-Riesz basis for $\mathcal{H}$.
  \end{remark}
  Next example justifies the above remark.
  \begin{exa}\label{eg4.3}
  	Let $\mathcal{H}$ be a separable Hilbert space with orthonormal basis $\{e_n\}_{n\in\mathbb{N}}$. For $m\in\mathbb{N}$, define $\Lambda_{m},\Omega_{m}:\mathcal{H}\rightarrow\mathbb{C}^4$ by
  	\begin{align*}
  		\Lambda_{m}(h)=\begin{cases}(\langle h,e_{2}\rangle,\langle h,e_{4}\rangle,\langle h,e_{1}\rangle,0) \quad  &\text{ if } m=1\\
  			(\langle h,e_{2}\rangle,\langle h,e_{4}\rangle,\langle h,e_{3}\rangle,0)\quad &\text{ if } m=2\\
  			(\langle h,e_{2}\rangle,\langle h,e_{4}\rangle,\langle h,e_{5}\rangle,\langle h,e_{6}\rangle)\quad  &\text{ if } m=3\\
  			(\langle h,e_{m+3}\rangle,0,0,0)\quad &\text{ if } m\geq 4\\
  		\end{cases}\\
  		\Omega_{m}(h)=\begin{cases}(\langle h,e_{1}\rangle,\langle h,e_{3}\rangle,\langle h,e_{2}\rangle,0) &\text{ if } m=1\\
  			(\langle h,e_{1}\rangle,\langle h,e_{3}\rangle,\langle h,e_{4}\rangle,0) &\text{ if } m=2\\
  			(\langle h,e_{1}\rangle,\langle h,e_{3}\rangle,\langle h,e_{5}\rangle,\langle h,e_{6}\rangle)& \text{ if } m=3\\
  			(\langle h,e_{m+3}\rangle,0,0,0) &\text{ if } m\geq 4.
  		\end{cases}
  	\end{align*}
  	Then $\{\Lambda_{m}\}_{m\in\mathbb{N}}$ and $\{\Omega_{m}\}_{m\in\mathbb{N}}$ are $g$-exact frames for $\mathcal{H}$. Let $\sigma\subseteq\mathbb{N}$ be arbitrary. Then
  	\begin{align*}
  		\|h\|^2&=\sum_{m\in\mathbb{N}}|\langle h,e_{m}\rangle|^2\\
  	&\leq\sum_{m\in\sigma}\|\Lambda_{m}h\|^2+\sum_{m\in\sigma^c}\|\Omega_{m}h\|^2\\
  		&\leq3\sum_{m\in\mathbb{N}}|\langle h,e_{m}\rangle|^2\\
  		&=3\|h\|^2,\quad \forall h\in\mathcal{H}.
  	\end{align*}
  	Thus $\{\Lambda_{m}\}_{m\in\mathbb{N}}$ and $\{\Omega_{m}\}_{m\in\mathbb{N}}$ are weaving $g$-frames for $\mathcal{H}$.
  	
  	For $\sigma=\{1,2\}$, $\{\Lambda_{m}\}_{m\in\sigma}\cup\{\Omega_{m}\}_{m\in\sigma^c}$ is a $g$-frame for $\mathcal{H}$. Since $\{\Lambda_{m}\}_{m\in\sigma\setminus\{2\}}\cup\{\Omega_{m}\}_{m\in\sigma^c}$ is also a $g$-frame for $\mathcal{H}$, therefore, $\{\Lambda_{m}\}_{m\in\sigma}\cup\{\Omega_{m}\}_{m\in\sigma^c}$ is not a $g$-exact frame and hence not a $g$-Riesz basis for $\mathcal{H}$.
  \end{exa}

\section*{Acknowledgments}
The authors thank Dr. Lalit Kumar Vashisht, University of Delhi, for his valuable comments and suggestions, which improved the presentation of the paper. Aniruddha Samanta thanks University Grants Commission(UGC)  for the financial support in the form of the Senior Research Fellowship (Ref.No:  19/06/2016(i)EU-V; Roll No. 423206).

	\bibliographystyle{amsplain}

\end{document}